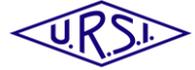

# Underwater Acoustic Channel Adaptive Estimation using $l_{21}$ Norms

Amer Aljanabi[1], Mohanad Abd Shehab[2], Osama Alluhaibi[3], Qasim Zeeshan Ahmed[1], and Pavlos Lazaridis[1]
1, Department of Engineering and Technology, University of Huddersfield, Huddersfield, HD1 3DH, UK
2, Electrical Engineering Department, Engineering College, Mustansiriyah University, Baghdad, IRAQ
3, School of Engineering, Electrical Engineering, University of Kirkuk, Kirkuk, IRAQ

**Abstract**

The problem of underwater acoustic (UWA) channel estimation is the non-uniform sparse representation that may increase the algorithm complexity and the required time. A mathematical framework utilizing $l_{21}$ constraint with two-dimensional frequency domain is employed to enhance the channel estimation. The frame work depends on both main and auxiliary channel information. The simulation results have been demonstrated that the proposed estimation method can improve some problems that are achieved with other norms like $l_1$. Furthermore, it can achieve a better performance in terms of mean square error (MSE) and execution time.

**Keywords:** Underwater acoustic channel, adaptive estimation, $l_1$ and $l_{21}$ norms.

## 1 Introduction

The primary purposes of underwater acoustic sensor networks include pollution monitoring, oceanographic data collection, disaster prevention, tactical surveillance, underwater exploration, and port security [1]. However, their drawbacks include limited battery life, difficult operating conditions (caused by a large number of multipaths), propagation delay spread due to the slow speed of sound, and susceptibility to Doppler shift between the source and destination. These have resulted in low data collection rates and limited distance communication [2]. Cooperative communication allows multiple sensors to form a distributed cooperative sensor network, enabling them to achieve spatial diversity which in turn helps save transmission power by combating the severe signal attenuation encountered over long distances. It can be difficult maintaining a lineof- sight (LOS) due to the movement of the platform caused by ambient disturbances or propulsion. There are therefore many potential applications involving mobile platforms, such as the one involving underwater robotics [3]. It also extends network coverage by providing LOS between the source and sensors, and between sensors and destination, which results in providing multiple communication links for higher data collection rates [4].

Two categories the current adaptive based channel estimation methods have mainly divided. Firstly, introducing a proportional step size into each filter coefficient the convergence rate by using the proportional adjustment, which uses the sparsity of the impulse response, such as, the improved lease mean square (LMS) algorithm [5-8]. Secondly, the sparse norm constraint, it is more much better and more precisely, a sparse penalty term ($l_{21}$ norm or $l_1$ norm) has been introduced into the LMS algorithm, in order to speeds up the convergence of
small coefficients [9, 10]. As a result, in the deep underwater environment the cluster-sparse distribution of the channel impulse response (CIR) cannot be fully exploited by the existing sparse adaptive algorithms. The best solution for this issue, we introduced a uniform $l_{21}$ norm constraint to the adaptive algorithms [11], the channel taps could be uniformly group without an overlapped, the algorithm apply $l_{21}$ norm constraint to the group and $l_1$ norm constraint between two groups [12]. The main issue is propagation models for UWA and tracking the UWA channel in underwater depths in real time. There are many complementary approaches proposed toward water acoustic channel estimation and the most importing thing is the real-time channel tracking remain a bottleneck [13], to solve these issues we should have a well knowledge about two underwater acoustic channel properties [14]. The long time-varying delay spread due to the moving ocean surface and primary and secondary multipath reflection [15]. The second reason is unpredictable high-energy due to oceanographic events [16]. In this work, many techniques are combined to support real-time channel estimation for intermediate distance and shallow water depths. These techniques include: managing non-uniform with sparse channel matrices, exploiting the main and auxiliary input data, and applying various regularization constraints to the optimization framework to provide good signal recovery. We will discourse more detail and challenges in different types of multipath arrivals as well [17,18].

## 2 Motivations and Challenges

The motivations of this paper revolve around providing channel state information (CSI). This is particularly challenging in underwater communications because the received underwater signal is often made up of many multipath components, each of which conveys very low energy. There are several complementary approaches have been discourse toward UWA channel tracking remain a bottleneck. To be more specific, these challenges
are addressed by two well-known properties of the underwater acoustic channel. The first one is long time-varying spread due to multipath reflection from moving ocean surface and static bottom as shown in figure (1).

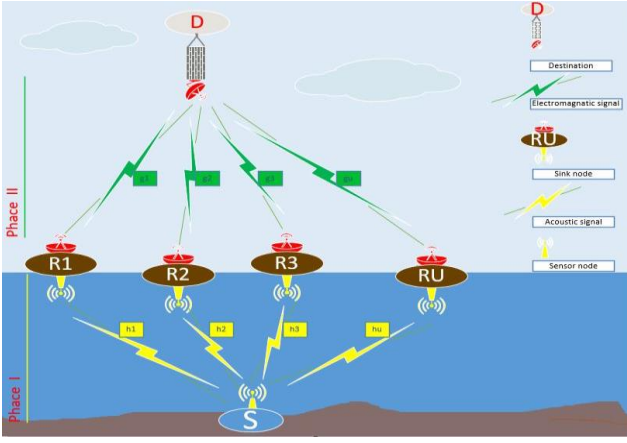

**Figure 1.** An underwater cooperative communication system

## 3 Channel Estimation Algorithm

The input data matrix U is a two-dimensional non-uniform matrix that can be bounded with L (number of the Doppler frequencies) × K (number of delay taps). Exploiting the 2D Fourier transform which can transform the channel estimation problem to spectral sampling problem as:

$$U = FH \quad (1)$$

where F is the 2D Fourier transform that can easily achieved with DFT transform and H it is the channel impulse response. In general, the received UWA signal will be disrupted with noise, i.e,

$$U_n = FH + N \quad (2)$$

Where Un denoted the noisy input data matrix and N is an additive noise.
Consider the auxiliary measured version of U as below

$$U_{aux} = RFH + N \quad (3)$$

where R is a binary random sub sampling matrix that permits for any positions of a random selection in the 2D Fourier domain. It can lead to sampling rate variations. The sampling rate S % can be formed depending on the dimensionality of the input matrix

$$S\% = \frac{KL \times S}{100} \quad (4)$$

As a result, to attain better generalization to the channel estimation, it can base on many regularizations that minimize the empirical error, i.e, it will be beneficial to take into account different penalties for different errors. These penalties are based on the Lp-norm like $l_1$, $l_2$ and $l_{21}$-norms as follow:

$$\|X\|_2 = \sqrt{\sum_{i=1}^{n}(X_i)^2}. \quad (5)$$

$$\|X\|_1 = \sum_{i=1}^{n}|X_i| \quad (6)$$

$$\|X\|_{21} = \sum_i (\sum_j X_{ij}^2)^{1/2} \quad (7)$$

where $l_2$ is easy to be achieved but it cannot minimize the error especially for the noisy data. $l_1$ can reject noisy data efficiently with the advantage of sparsity but with complex convex solution while $l_{21}$ is often applied to tackle and overtake the difficulty of data. Finally, the resultant channel objective functions can be modelled with two main forms as in equations (8) and (9):

$$\min_{H} \|U_{aux} - RFH\|_2^2 \quad s.t \quad \|H\|_1 \leq \sigma \quad (8)$$

$$\min_{H} \|JH\|_{21} \quad s.t \quad \|U_n - F.H\|_2^2 \leq \sigma \quad (9)$$

Where $\sigma$ = stander deviation of the noise

## 4 Numerical Result

In this section, the channel estimation performance of the proposed non-uniform $l_{21}$ is compared with that of the $l_1$ for a simple channel model that includes the following parameters:
L= number of the Doppler frequencies = 11 × K= number of delay taps = 200, frequency resolution = 25Hz, time resolution = 0.05 sec, H = [361 × 500].
Figures (2) and (3) show the MMSE for different window length and sampling ratios. It is clear that algorithm that based on $l_{21}$ is of better response than that of $l_1$. Also, the time required for achieving $l_{21}$ is two-third of that required for $l_1$.

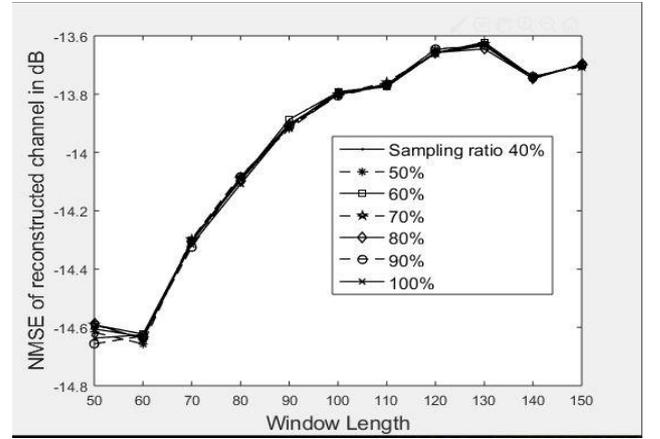

**Figure 2.** Channel estimation MSE at SNR = 10 dB with norm $l_1$

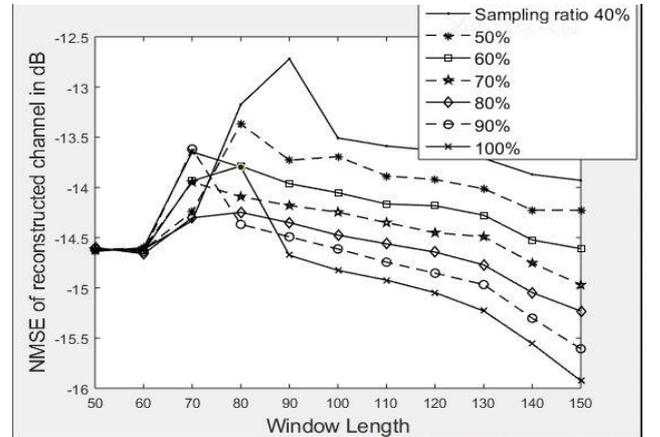

**Figure 3.** Channel estimation MSE at SNR = 10 dB with norm $l_{21}$

## 5 Conclusion

The underwater acoustic channel estimation accuracy and speed are vital in communication system. However, it is clear that the responsibility for improving underwater acoustic communication lies with channel estimation and be achieved via adopting robust algorithms. Herein, we are proposed different objective functions with norm 1 and norm 21 penalties that can minimize the MMSE error and increase the estimation speed. The results for some cases approve that the norm 21 is better than norm 1 in both MMSE and execution time.